\documentclass[12pt,letterpaper]{amsart}

\usepackage{amssymb,amscd}
\usepackage[mathscr]{eucal}

\numberwithin{equation}{section}

\newtheorem{theorem}{Theorem}
\newtheorem{proposition}{Proposition}
\newtheorem{lemma}{Lemma}
\newtheorem{definition}{Definition}
\numberwithin{theorem}{section}
\numberwithin{proposition}{section}
\numberwithin{lemma}{section}
\numberwithin{definition}{section}

\newcommand{\C}{\mathbb{C}}
\newcommand{\Z}{\mathbb{Z}}

\renewcommand{\L}{\mathcal{L}}
\newcommand{\CO}{\mathcal{O}}
\newcommand{\CP}{\mathcal{P}}
\newcommand{\CQ}{\mathcal{Q}}
\newcommand{\CX}{\mathsf{X}}
\newcommand{\CH}{\mathcal{H}}
\newcommand{\CC}{\mathcal{C}}
\newcommand{\CG}{\mathcal{G}}
\newcommand{\QQ}{\mathsf{Q}}
\newcommand{\LL}{\mathbb{L}}
\newcommand{\DD}{\mathcal{D}}
\newcommand{\bull}{\bullet}
\renewcommand{\o}{\otimes}
\newcommand{\Dual}{*}

\DeclareMathOperator{\End}{End}
\DeclareMathOperator{\iso}{\mathfrak{iso}}
\DeclareMathOperator{\so}{\mathfrak{so}}
\newcommand{\h}{\mathfrak{h}}
\newcommand{\g}{\mathfrak{g}}
\DeclareMathOperator{\ad}{ad}
\DeclareMathOperator{\Ad}{Ad}
\newcommand{\half}{\tfrac{1}{2}}
\newcommand{\p}{\partial}
\renewcommand{\*}{\cdot}

\renewcommand{\]}{{]\!]}}
\renewcommand{\[}{{[\![}}
\newcommand{\m}{\mathsf{m}}
\DeclareMathOperator{\Sl}{\mathfrak{sl}}
\DeclareMathOperator{\MC}{\mathcal{MC}}
\DeclareMathOperator{\gr}{gr}
\newcommand{\tint}{{\textstyle\int}}
\newcommand{\Om}{\Omega}
\newcommand{\om}{\omega}
\newcommand{\eps}{\varepsilon}
\newcommand{\RR}{\mathsf{T}}
\newcommand{\Q}{\mathbb{Q}}

\begin{document}

\title[A Darboux theorem for Hamiltonian operators]
{A Darboux theorem for Hamiltonian operators in the formal calculus
of variations}

\author[E. Getzler]{Ezra Getzler}

\address{RIMS, Kyoto University, Kitashirakawa Oiwake-cho, Sakyo-ku, Kyoto,
Japan}

\address{Department of Mathematics, Northwestern University, Evanston,
IL 60208, USA}

\email{getzler@math.northwestern.edu}

\dedicatory{This paper is dedicated to Roger Richardson \emph{in memoriam}.}

\maketitle

\thispagestyle{empty}

The Darboux theorem states that all symplectic structures on an affine
space are isomorphic. In the formal calculus of variations, symplectic forms
are replaced by Hamiltonian operators, which are systems of ordinary
differential equations satisfying a complicated quadratic constraint. It is
natural to ask whether in this setting an analogue of the Darboux theorem
holds.

The problem is considerably simplified if one restricts attention to formal
deformations of a given Hamiltonian operator $H$. The study of the moduli
space of deformations is then controlled by a differential graded (dg) Lie
algebra, the Schou\-ten Lie algebra, with differential induced by $H$. The
problem of calculating the cohomology of this dg Lie algebra was posed by
Olver \cite{Olver}.

Let $V$ be a finite-dimensional vector space with basis $e_a$, and let
$\eta=\eta^{ab}e_a \o e_b$ be a non-degenerate bilinear form on
$V^\Dual$. In this paper, we answer Olver's question for the hydrodynamic
Hamiltonian operator $H=\eta^{ab}\p$. We show that the associated dg Lie
algebra is formal; that is, it is quasi-isomorphic to its cohomology.

As a consequence, the set of possible normal forms of a deformation of the
Hamiltonian operator $\eta^{ab}\p$ is very easy to calculate: the answer is
related to a result of Dubrovin and Novikov \cite{DN}, who showed that an
operator of the form
$$
\eta^{ab} \p + A^{ab}_c t^c + B^{ab}
$$
is Hamiltonian if and only if $[e^a,e^b]=A^{ab}_ce^c$ is a Lie bracket on
$V^\Dual$ with the metric $\eta$ as a Killing form, and $B^{ab}e_a\wedge
e_b$ is a two-cocycle for this Lie bracket.

We actually work in a global setting, replacing the vector space $V$ by a
complex manifold $X$ with flat contravariant metric $\eta$. (A
contravariant metric is a non-degenerate symmetric bilinear form on the
cotangent bundle; we do not, indeed cannot, impose any positivity condition
on $\eta$, since we work with complex manifolds and holomorphic sections.)

We utilize three main concepts in this paper:
\begin{enumerate}
\item dg Lie algebras concentrated in degrees $[-1,\infty)$, such as the
Schouten algebra --- these give a convenient language for describing
deformation problems;
\item the Deligne 2-groupoid associated to such a dg Lie algebra, which
represents the moduli of formal deformations;
\item a refined version of the Schouten bracket in the formal calculus
of variations, due to V.~O.~Soloviev \cite{Soloviev}.
\end{enumerate}
These three ingredients are explained respectively in Sections 1, 2 and
4. In Section~3, we discuss the structure of the Deligne 2-groupoids which
arise in the deformation theory of hydrodynamic Hamiltonian operators. Our
main result is stated in Section 5, and proved in Section 6.

Throughout this paper, the summation convention is understood: indices
$a,b,\dots$ once as superscript and once as subscript in a formula are to
be summed over.

\subsection*{Acknowledgements}
This paper contributes to an area of mathematics in which my teacher Roger
Richardson was a pioneer \cite{NR}. He first introduced me to the beautiful
applications of algebra in differential geometry.

Youjin Zhang brought the problem of proving a Darboux theorem for
Hamiltonian operators to my attention. I thank him, B.~Dubrovin, E.~Frenkel
and P.~Olver for useful discussions.

The ``Geometry of String Theory'' project of RIMS provided support while
this paper was written; I am especially grateful to K. Saito for helping
make my visit to RIMS so enjoyable.

Additional support was provided by the NSF through grant DMS-9704320.

\section{Poisson tensors on supermanifolds and the Schouten bracket}
\label{poisson}

In this section, we recall the elements of the theory of Poisson
supermanifolds. This theory differs a little from that of Poisson
manifolds, since the Poisson tensor on a supermanifold may have either even
or odd parity.

\subsection{Poisson tensors on supermanifolds}
Let $\C^{m|n}$ be the superspace with $m$ even and $n$ odd coordinates; if
$U$ is an open subset of $\C^{m|n}$ and $|U|$ is the underlying open subset
of $\C^m$, we have $\CO(U) \cong \CO(|U|) \o \Lambda(\C^n)^\Dual$.
\begin{definition}
A \emph{$\nu$-Poisson tensor} on $U$ is a two-tensor
$$
P = P^{ab}\,\p_a\o\p_b \in \CO(U) \o (\C^{m|n}){}^{\o2} ,
$$
of total degree $\nu\in\Z/2$ (i.e. $|P^{ab}|=|a|+|b|+\nu$), such that
$P^{ba}+(-1)^{|a||b|+\nu} P^{ab}=0$ and
\begin{equation} \label{Poisson}
\sum_{\substack{\textup{cycles in} \\ (b,c,d)}} (-1)^{|b|(|a|+|c|+\nu)}
P^{ba} \, \p_aP^{cd} = 0 .
\end{equation}
\end{definition}

A $\nu$-Poisson tensor $P$ on $U$ defines a Poisson bracket on $\CO(U)$, by
the formula
$$
\{u,v\} = - (-1)^{(|b|+\nu)|u|} P^{ab} \, \p_au \, \p_bv .
$$
The symmetry of $P$ is equivalent to skew-symmetry of the bracket,
$$
\{u,v\} + (-1)^{(|u|+\nu)(|v|+\nu)} \{v,u\} = 0 ,
$$
while \eqref{Poisson} is equivalent to the Jacobi rule,
$$
\{u,\{v,w\}\} - (-1)^{(|u|+\nu)(|v|+\nu)} \{v,\{u,w\}\} = \{\{u,v\},w\} .
$$
We conclude that the space of $\nu$-Poisson tensors is invariant under
change of coordinates; thus, we may define a $\nu$-Poisson tensor on a
complex supermanifold as tensor which is a $\nu$-Poisson tensor for some
atlas.
\begin{definition}
A holomorphic \emph{$\nu$-Poisson} supermanifold $(M,P)$ is a complex
supermanifold $M$ together with a $\nu$-Poisson tensor $P$ on $M$. If the
Poisson tensor is non-degenerate, we call $(M,P)$ a holomorphic
\emph{$\nu$-symplectic} supermanifold.
\end{definition}

\subsection{The Schouten bracket and supermanifolds}
If $X$ is a manifold, let $\Om X$ be the $1$-symplectic supermanifold
obtained by forming the cotangent bundle $T^\Dual\!X$ over $X$, which is a
symplectic manifold, applying the functor $\Pi$ which reverses the parity
of the fibres, and taking the underlying supermanifold. Let $\pi:\Om X\to
X$ be the projection.

Let $t^a$, $1\le a\le m$, be coordinates on an open subset of $X$, and let
$\theta_a$ be the dual coordinates along the fibres of $\Om X$; let
$\p_a=\p/\p t^a$ and $\p^a=\p/\p\theta_a$ be the corresponding vector
fields. The Poisson tensor (or more accurately, $1$-Poisson tensor) of $\Om
X$ is given by the formula
\begin{equation} \label{Om}
P = \p^a \o \p_a + \p_a \o \p^a .
\end{equation}
The sheaf $\pi_*\CO_{\Om X}$ is isomorphic to the graded sheaf
$\Lambda=\Lambda T_X$ of multivectors on $X$, and this isomorphism
identifies the Poisson bracket on $\Om X$ with the Schouten bracket $[-,-]$
on $\Lambda$.

The Poisson bracket on the $\Z/2$-graded sheaf $\CO_{\Om X}$ has odd
degree, while we prefer to work with a $\Z$-grading on $\Lambda$ such that
the Schouten bracket has degree $0$. To this end, we define the degree $p$
summand $\Lambda^p$ of $\Lambda$ to be $\Lambda^{p+1}T_X$. Taking this
shift of degree into account, the formula for the Schouten bracket becomes
$$
[u,v] = (-1)^{|u|} \p_au \, \p^av - \p^au \, \p_av .
$$
If $Q\in\Gamma(X,\Lambda^1)=\Gamma(X,\Lambda^2T_X)$, define an operation
$\{u,v\}_Q$ on $\CO$ by the formula
$$
\{u,v\}_Q = [[Q,u],v] .
$$
In local coordinates $Q=\half Q^{ab}\theta_a\,\theta_b$, we have
$$
[Q,u] = Q^{ab} \theta_a \, \p_bu ,
$$
and hence $\{u,v\}_Q=Q^{ab} \p_au \, \p_bu$. If $Q$ is a Poisson tensor,
this is the Poisson bracket associated to $Q$.
\begin{proposition} \label{L}
The following conditions on a section $Q$ of $\Lambda^1=\Lambda^2T_X$ are
equivalent:
\begin{enumerate}
\item $Q$ is a Poisson tensor on $X$;
\item $[Q,Q]=0$;
\item the operation $\delta_Q=[Q,-]$ is a differential on the sheaf of
graded Lie algebras $\Lambda$;
\item the operation $\{u,v\}_Q=[\delta_Qu,v]$ on $\CO$ is a Lie bracket.
\end{enumerate}
\end{proposition}
\begin{proof}
In local coordinates, the formula $[Q,Q]=0$ becomes Eq.\ \eqref{Poisson}
for the tensor $Q$ on $X$; thus (1) and (2) are equivalent.

The Jacobi rule for graded Lie algebras shows that
$$
\delta_Q\delta_Qa = [Q,[Q,a]] = \half [[Q,Q],a] .
$$
Thus $\delta_Q$ is a differential on $\Lambda$ if and only if $[Q,Q]=0$.

The bracket $\{u,v\}_Q$ is skew-symmetric:
$$
\{u,v\}_Q + \{v,u\}_Q = [ \delta_Qu,v ] + [ \delta_Qv,u ] = \delta_Q[u,v] =
0 .
$$
As for the Jacobi rule, we have
\begin{align*}
\{u,\{v,w\}\} - \{v,\{u,w\}\} &= [\delta_Qu,[\delta_Qv,w]] -
[\delta_Qv,[\delta_Qu,w]] \\
&= [[\delta_Qu,\delta_Qv],w] = \{\{u,v\},w\} - [[\delta_Q\delta_Qu,v],w] .
\end{align*}
The anomalous term $-\half[[\delta_Q\delta_Qu,v],w]$ vanishes for all $u$,
$v$ and $w\in\Gamma(U,\CO)$ if and only if $\delta_Q$ is a differential.
\end{proof}

Let $(X,Q)$ be a Poisson manifold. We denote the sheaf of dg Lie algebras
$\Lambda$, with differential $\delta_Q$, by $\Lambda_Q$. For example, if
$(X,Q)$ is a symplectic manifold, then the complex of sheaves underlying
$\Lambda_Q$ is isomorphic to the de Rham complex, and the Poisson
cohomology is isomorphic to the trivial sheaf $\C$, with vanishing Lie
bracket.

\section{Graded Lie algebras and the Deligne 2-groupoid}

Goldman and Millson \cite{GM} have developed an approach to deformation
theory based on a functor from nilpotent dg Lie algebras concentrated in
degrees $[0,\infty)$ to groupoids $\CC(\g)$, called the Deligne groupoid.

The dg Lie algebra controlling the deformation theory of Poisson brackets
is the Schouten Lie algebra, which is concentrated in degrees
$[-1,\infty)$; thus, the theory of the Deligne groupoid does not apply. It
turns out that the deformation theory is best understood by means of a
2-groupoid, whose definition generalizes that of the Deligne
groupoid.\footnote{We have learned that this 2-groupoid was proposed by
Deligne in a letter to Breen (February, 1994); it is also alluded to in
Section 3.3 of Kontsevich \cite{Kontsevich}.}

In this section, all dg Lie algebras $\g$ are concentrated in degrees
$[-1,\infty)$.

\subsection{The Deligne groupoid}

We now recall the definition of the Deligne groupoid. There is a sequence
of elements $F_n(x,y)$ of degree $n$ in the free Lie algebra on two
generators $x$ and $y$ such that if $X$ and $Y$ are elements of a nilpotent
Lie algebra $\g$ of $N$ steps, we have
$$
\exp(X)\exp(Y) = \exp\Bigl( \sum_{n=1}^N F_n(X,Y) \Bigr)
$$
in the associated simply connected Lie group $G$; for example,
$F_1(X,Y)=X+Y$ and $F_2(X,Y)=\half[X,Y]$. We may identify the Lie group $G$
with the manifold $\g$ with deformed product
$$
X \* Y = \sum_{n=1}^N F_n(X,Y) .
$$
Denote the resulting functor from nilpotent Lie algebras to Lie groups by
$\exp(\g)$.

\begin{definition}
If $\g$ is a dg Lie algebra, the set $\MC(\g)$ of \emph{Maurer-Cartan}
elements of $\g$ is the inverse image $Q^{-1}(0)\subset\g^1$ of the
quadratic map $Q:\g^1\to\g^2$ defined by the formula $Q(A) =
dA+\half[A,A]$.
\end{definition}

Thus, $A$ is a Maurer-Cartan element if and only the operator $d_Au = du +
[A,u]$ is a differential on $\g$.

The subspace $\g^0$ of $\g$ is a nilpotent Lie algebra, and the group
$\exp(\g^0)$ acts on $\g^1$ by the formula
\begin{equation} \label{action}
\exp(X) \* A = A - \sum_{n=0}^\infty \frac{\ad(X)^n}{(n+1)!} \, d_AX ;
\end{equation}
this is the affine action which corresponds to gauge transformations
$$
d_{\exp(X)\*A}=\Ad(\exp(X))d_A .
$$
Since $Q(\exp(X)\*A)=Q(A)$, this action preserves the subset
$\MC(\g)\subset\g^1$.
\begin{definition}
The \emph{Deligne groupoid} $\CC(\g)$ of $\g$ is the groupoid associated
to the group action $\exp(\g^0)\times\MC(\g)\to\MC(\g)$.
\end{definition}

The sets of objects and morphisms of the Deligne groupoid are $\MC(\g)$ and
$\exp(\g^0)\times\MC(\g)$; its source and target maps are $s(\exp(X),A)=A$,
and $t(\exp(X),A)=\exp(X)\*A$, its identity is $A\mapsto(\exp(0),A)$, and
its composition is
$$
(\exp(Y),\exp(X)\*A) \* (\exp(X),A) = (\exp(Y)\exp(X),A) .
$$
The Deligne groupoid is a natural generalization of Lie's correspondence
$\exp$ between nilpotent Lie algebras and simply connected nilpotent Lie
groups, which is the case where $\g$ is concentrated in degree $0$.

Even if $\g$ is not nilpotent, we may consider the Deligne groupoid with
coefficients in $\g\o\m$, where $\m$ is a nilpotent commutative algebra.

If $\CG$ is a groupoid, let $\pi_0(\CG)$ be the set obtained by quotienting
of the set of objects of $\CG$ by the equivalence relation $x\sim y$
whenever there is a morphism between $x$ and $y$. If $\g$ is a nilpotent dg
Lie algebra, we will write $\pi_0(\g)$ for $\pi_0(\CC(\g))$. Much of
deformation theory may be reformulated as the study of the sets
$\pi_0(\g\o\m)$.

The following result is proved by exactly the same method as Theorem 2.4 of
Goldman and Milsson \cite{GM}.
\begin{theorem} \label{gm}
Let $\h=F^1\h\supset F^2\h\supset\dots$ and
$\tilde{\h}=F^1\tilde{\h}\supset F^2\tilde{\h}\supset\dots$ be filtered dg
Lie algebras (that is, $dF^i\h\subset F^i\h$ and $[F^i\h,F^j\h]\subset
F^{i+j}\h$, and likewise for $\tilde{\h}$) such that $F^N\h$ and
$F^N\tilde{\h}$ vanish for sufficiently large $N$, and let
$f:\h\to\tilde{\h}$ be a morphism of filtered dg Lie algebras which induces
weak equivalences of the associated chain complexes
$$
\gr^i f:F^i\h/F^{i+1}\h \longrightarrow
F^i\tilde{\h}/F^{i+1}\tilde{\h} .
$$
Then $f$ induces a bijection $\pi_0(f) : \pi_0(\h) \to \pi_0(\tilde{\h})$.
\end{theorem}

\subsection{2-groupoids}

The category of groupoids is a monoidal category, where $\CG\o\CH$ is the
product $\CG\times\CH$ of the groupoids $\CG$ and $\CH$.
\begin{definition}
A \emph{2-groupoid} is a groupoid enriched over the monoidal category of
group\-oids.
\end{definition}

We see that a 2-groupoid $\CG$ has a set $G_0$ of objects, and for each
pair of objects $x,y\in G_0$, a groupoid of morphisms $\CG(x,y)$, and that
there are product maps
\begin{equation} \label{Horizontal}
\CG(x,y) \times \CG(y,x) \longrightarrow \CG(x,z) ,
\end{equation}
satisfying the usual conditions of associativity for a category.

The 2-morphisms of a 2-groupoid are the morphisms of the groupoids
$\CG(x,y)$. There are two compositions defined on the 2-morphisms: the
\emph{horizontal} composition of \eqref{Horizontal} and the \emph{vertical}
composition, which is composition inside the groupoid $\CG(x,y)$.

If $\CG$ is a 2-groupoid, let $\pi_1(\CG)$ be the groupoid whose objects
are those of $\CG$, and such that the set of morphisms $\pi_1(\CG)(x,y)$
equals $\pi_0(\CG(x,y))$. Let $\pi_0(\CG)$ equal $\pi_0(\pi_1(\CG))$.

If $x$ is an object of $\CG$, let $\pi_1(\CG,x)$ be the automorphism group
$\pi_1(\CG)(x,x)$, and let $\pi_2(\CG,x)$ be the automorphism group of the
identity of $x$ in the groupoid $\CG(x,x)$. The group $\pi_2(\CG,x)$ is
abelian, for the same reason as $\pi_2(X,x)$ is abelian for a topological
space $X$: it carries two products, horizontal and vertical, satisfying
$(a\circ_hb)\circ_v(c\circ_hd)=(a\circ_vc)\circ_h(b\circ_vd)$.

\begin{definition}
A \emph{weak equivalence} $\varphi:\CG\to\CH$ of 2-groupoids is a
homomorphism such that $\pi_0(\varphi)$ is an isomorphism of sets and for
each object $x\in G_0$,
$$
\pi_i(\varphi,x) : \pi_i(\CG,x)\to\pi_i(\CH,\varphi(x))
$$
is an isomorphism of groups for all $x\in G_0$ and $i=1,2$.
\end{definition}

With this notion of weak equivalence and suitable notions of cofibration
and fibration, the category of 2-groupoids is a closed model category
(Moerdijk and Svensson \cite{MS}).

\subsection{The Deligne 2-groupoid}

We now show how the Deligne groupoid of a nilpotent dg Lie algebra $\g$ is
the underlying groupoid of a 2-groupoid, which we denote by $\CC(\g)$; if
$\g$ happens to vanish in degree $-1$, this 2-groupoid is identical to the
Deligne groupoid of $\CG$. (Thus, the use of the same notation for the
Deligne 2-groupoid and Deligne groupoid should cause no difficulty.)

Given an element $A\in\g^1$, define a bracket $\{u,v\}_A$ on $\g^{-1}$, by
the formula
\begin{equation} \label{gA}
\{u,v\}_A = [d_Au,v] .
\end{equation}
The proof of the following proposition is the same as the proof of the
equivalence of conditions (3) and (4) in Proposition \ref{L}.
\begin{proposition} \label{l}
The bracket $\{u,v\}_A$ makes $\g^{-1}$ into a Lie algebra if and only if
$A\in\MC(\g)$.
\end{proposition}

If $A\in\MC(\g)$, we denote the Lie algebra $\g^{-1}$ with bracket
$\{u,v\}_A$ by $\g_A$. The nilpotence of $\g$ implies that that $\g_A$ is
nilpotent. If $u\in\g_A\cong\g^{-1}$, denote the corresponding element of
the group $\exp(\g_A)$ by $\exp_A(u)$.

Since $d_A\{u,v\}_A = [d_Au,d_Av]$, the linear map $d_A:\g_A\to\g^0$ is a
morphism of Lie algebras. Thus, the group $\exp(\g_A)$ acts on $\exp(\g^0)$
by right translation:
$$
\exp_A(u) \* \exp(X) = \exp(X)\exp(d_Au) .
$$
Given a pair $A,B$ of elements of $\MC(\g)$, define $\CC(\g)(A,B)$ to be
the groupoid associated to this group action. The set of 2-morphisms of
$\CC(\g,\m)$ may be identified with $\g^{-1}\times\g^0\times\MC(\g)$; we
denote its elements by $(\exp_A(u),\exp(X),A)$, where $u\in\g^{-1}$,
$X\in\g^0$ and $A\in\MC(\g)$. The internal (or vertical) composition of
2-morphisms is given by the formula
$$
\bigl( \exp_A(v) , \exp(X)\exp(d_Au) , A \bigr) \circ_v \bigl( \exp_A(u) ,
\exp(X) , A \bigr) = \bigl( \exp_A(u)\exp_A(v) , \exp(X) , A \bigr) .
$$

To complete the definition of the Deligne 2-groupoid, it remains to define
the horizontal composition
\begin{equation} \label{horizontal}
\CC(\g)(B,C) \times \CC(\g)(A,B) \longrightarrow \CC(\g)(A,C) .
\end{equation}
Given $X\in\g^0$ and $A\in\MC(\g)$, there is an isomorphism of Lie algebras
$$
e^{\ad(X)} : \g_{\exp(X)\*A} \longrightarrow \g_A ,
$$
with inverse $e^{-\ad(X)}$. Suppose $A\in\MC(\g)$ and $X,Y\in\g^0$, with
$B=\exp(X)\*A$ and $C=\exp(Y)\*B$. Then the horizontal composition of
2-mophisms is given by the formula
\begin{multline*}
\bigl( \exp_{\exp(X)\*A}(v) , \exp(Y) , \exp(X)\*A \bigr) \circ_h \bigl(
\exp_A(u) , \exp(X) , A \bigr) \\ = \bigl(
\exp_A\bigl(e^{-\ad(X)}v\bigr)\exp_A(u) , \exp(Y)\exp(X) , A \bigr) .
\end{multline*}

If $\g$ is a nilpotent dg Lie algebra and $A\in\MC(\g)$, we will write
$\pi_1(\g,A)$ and $\pi_2(\g,A)$ for $\pi_1(\CC(\g),A)$ and
$\pi_2(\CC(\g),A)$.

\begin{theorem} \label{GM}
Let $\g$ and $\tilde{\g}$ be dg Lie algebras concentrated in degrees
$[-1,\infty)$, and let $\m$ be a nilpotent commutative algebra. A weak
equivalence $f:\g\to\tilde{\g}$ of dg Lie algebras induces a weak
equivalence of 2-groupoids $\CC(f\o\m) : \CC(\g\o\m) \to
\CC(\tilde{\g}\o\m)$.
\end{theorem}
\begin{proof}
By Theorem \ref{gm}, a weak equivalence $f:\g\to\tilde{\g}$ of dg Lie
algebras induces a bijection $\pi_0(f\o\m) : \pi_0(\g\o\m) \to
\pi_0(\tilde{\g}\o\m)$; indeed, $\g\o\m$ is filtered by subspaces
$F^i\g\o\m=\g\o\m^i$, and similarly for $\tilde{\g}$. It remains to prove
that $f$ induces bijections $\pi_i(\g\o\m,A)\cong\pi_i(\tilde{\g}\o\m,f(A)$
for all $A\in\MC(\g\o\m)$ and $i\in\{1,2\}$.

Given $A\in\MC(\g\o\m)$, define a dg Lie algebra
$$
\Om_A(\g\o\m) = \bigl( 0 \longrightarrow (\g\o\m)_A \xrightarrow{\ d_A\ }
\ker(d_A|\g^0\o\m) \longrightarrow 0 \bigr) ,
$$
where $(\g\o\m)_A$ is placed in degree $0$. The construction
$\Om_A(\g\o\m)$ behaves like a based loop space of $\g\o\m$ at $A$, in the
sense that
\begin{equation} \label{loop}
\CC(\g\o\m)(A,A)\cong\CC(\Om_A(\g\o\m))
\end{equation}
To prove this, we must first show that these groupoids have the same
objects, that is, that $\exp(X)\*A=A$ if and only if $d_AX=0$. If
$\exp(X)\*A=A$, we see from \eqref{action} that
$$
d_AX = - \sum_{n=1}^\infty \frac{\ad(X)^n}{(n+1)!} \, d_AX .
$$
It follows, by induction on $n$, that $d_AX\in\g^1\o\m^n$ for all $n>0$,
hence $d_AX=0$.  The remainder of the proof of \eqref{loop} is
straightforward.

That $\pi_1(f\o\m,A):\pi_1(\g\o\m,A)\to\pi_1(\tilde{\g}\o\m,f(A))$ is a
bijection now follows on applying Theorem \ref{gm} to the weak equivalence
of filtered dg Lie algebras
$$
\Om_A(f\o\m) : \Om_A(\g\o\m) \to \Om_{f(A)}(\tilde{\g}\o\m) .
$$

Finally, $\pi_2(f\o\m,A) : \pi_2(\g\o\m,A) \longrightarrow
\pi_2(\tilde{\g}\o\m,f(A))$ is a bijection, since $\pi_2(\g\o\m,A)\cong
H^{-1}(\g\o\m,d_A)$.
\end{proof}

If $\g$ is a dg Lie algebra, its cohomology $H(\g)$ is a dg Lie algebra
with vanishing differential.
\begin{definition}
A dg Lie algebra $\g$ is \emph{formal} if there exists a dg Lie algebra
$\tilde{\g}$ and weak equivalences of dg Lie algebras $\tilde{\g}\to\g$ and
$\tilde{\g}\to H(\g)$.
\end{definition}

If $\g$ is formal, Theorem \ref{GM} implies that the 2-groupoids
$\CC(\g,\m)$ and $\CC(H(\g),\m)$ are equivalent, and hence that the
2-groupoid $\CC(H(\g),\m)$ parametrizes normal forms for deformations of
the differential on $\g$. This motivates the following.
\begin{definition} \label{Darboux}
A deformation problem is \emph{Darboux} if it is controlled by a formal
dg Lie algebra $\g$.
\end{definition}

\section{Examples of Deligne 2-groupoids}

We now illustrate the Deligne 2-groupoid in two examples: the deformation
theory of Poisson tensors, and a graded Lie algebra which occurs in the
deformation theory of Hamiltonian operators of hydrodynamic type.

\subsection{Deformation of Poisson tensors}
Let $(X,Q)$ be a Poisson manifold, and let $\g$ be the dg Lie algebra
$\g=\Gamma(X,\Lambda)$ with differential $d_Q$.

Given an integer $n$, let $\m_n$ be the nilpotent ring
$\hbar\,\C[\hbar]/(\hbar^{n+1})$. The Maurer-Cartan elements of $\g\o\m_n$
are the $n$th order deformations
$$
\QQ = Q + \sum_{k=1}^n \hbar^k Q_k + O(\hbar^{n+1}) , \quad Q_k \in \g^1 ,
$$
of the Poisson tensor such that $[\QQ,\QQ]=O(\hbar^{n+1})$. The Lie
algebra $\g^0\o\m_n$ may be identified with the Lie algebra of formal
vector fields
$$
\CX = \sum_{k=1}^n \hbar^k X_k + O(\hbar^{n+1}) , \quad X_k \in
\Gamma(X,\Lambda^1) ,
$$
and $\exp(\g^0\o\m_n)$ with the group of formal diffeomorphisms; thus
$\pi_0(\g\o\m_n)$ is the set of equivalence classes of $n$th order
deformations $\QQ$ of the Poisson bracket $Q$ modulo formal
diffeomorphisms. For $i=1,2$, we have
$$
\pi_i(\g\o\m_n,\QQ) \cong \exp(H^{1-i}(\Gamma(X,\Lambda)\o\m_n,\delta_\QQ))
;
$$
in particular, $\pi_2(\g\o\m_n,\QQ)$ is the space of Casimirs of $\QQ$.

The deformation theory of an affine symplectic manifold $(V,Q)$ is Darboux
in the sense of Definition \ref{Darboux}: its controlling dg Lie algebra
$\Gamma(V,\Lambda_Q)$ has cohomology $\C[1]$, and hence is formal; in this
way, we recover a formal version of the usual Darboux theorem. From this
example, we see how powerful formality is: it allows the calculation of the
homotopy type of the Deligne 2-groupoid (in this case, $K(\C,2)$) in a
straightforward way.

\subsection{A Deligne 2-groupoid associated to a Euclidean vector space}
\label{hydro-example}

We now consider the Deligne 2-groupoids of a class of graded Lie algebras
associated to Euclidean vector spaces.

If $(V,\eta)$ is a Euclidean vector space, the odd superspace $\Pi V^\Dual$
is symplectic (i.e.\ $0$-sym\-plec\-tic). If $t^a$ is a coordinate system on
$V$ (that is, a basis of $V^\Dual$) and $\theta_a$ is the dual coordinate
system on $\Pi V^\Dual$, the symplectic form on $\Pi V^\Dual$ equals
$$
\om = \eta^{ab} d\theta_a \, d\theta_b .
$$

The Lie algebra $\h(V,\eta)$ of Hamiltonian vector fields on $\Pi V^\Dual$
is a $\Z$-graded Lie algebra: the Poisson bracket has degree $-2$ (with
respect to the degree in the generators $\theta_a$ of $\CO_{\Pi
V^\Dual}=\Lambda V$), so the $\Z$-grading is defined by assigning to a
Hamiltonian vector field its degree of homogeneity minus $1$. Equivalently,
this equals the degree of homogeneity of the corresponding Hamiltonian
minus $2$; thus
$$
\h^p(V,\eta) \cong \begin{cases} \Lambda^{p+2}V & p\ge-1 , \\ 0 & p<-1
. \end{cases}
$$
Using the Hamiltonians to represent the corresponding Hamiltonian vector
fields, the bracket of elements $\alpha\in\h^p(V,\eta)$ and
$\beta\in\h^q(V,\eta)$ is
$$
\{ \alpha , \beta \} = (-1)^{p+1} \, \eta_{ab} \, \p^a \alpha \, \p^b \beta
.
$$

The graded vector space $\CO[1]$ is a graded module for the graded Lie
algebra $\h(V,\eta)$, with $\CO[1]^p\cong\Lambda^{p+1}V$; the action of
$\alpha\in\h^p(V,\eta)$ on $\tilde\beta\in\CO[1]$ is given by the formula
$$
\alpha \* \beta = - \eta_{ab} \, \p^a \alpha \, \p^b \tilde\beta .
$$
The sign is explained by the fact that we consider the module $\CO[1]$ and
not $\CO$.

Let $\g(V,\eta)=\CO[1]\rtimes\h(V,\eta)$ be the semidirect product of
$\h(V,\eta)$ with the abelian graded Lie algebra $\CO[1]$; thus,
$\g^p(V,\eta)$ is isomorphic to $\Lambda^{p+1}V\oplus\Lambda^{p+2}V$. We
denote elements of $\g^p(V,\eta)$ by $(\tilde\alpha,\alpha)$, where
$\tilde\alpha\in\CO[1]^p\cong\Lambda^{p+1}V$ and
$\alpha\in\h^p(V,\eta)\cong\Lambda^{p+2}V$.

The Lie subalgebra $\g^0(V,\eta)\subset\g(V,\eta)$ is isomorphic to
$\iso(V,\eta)$, the Lie algebra of infinitesimal Euclidean transformations
of $V$; $\g(V,\eta)$ is an analogue of $\iso(V,\eta)$ in the graded world.

Let $\m$ be a nilpotent commutative algebra. An element
$(\tilde\alpha,\alpha)$ of
$$
\g^1(V,\eta)\o\m\cong(\Lambda^3V\oplus\Lambda^2V)\o\m
$$
gives rise to a skew-symmetric operation on $(\C\oplus V^\Dual)\o\m$ by the
formula
$$
[(a,x),(b,y)]_{(\tilde\alpha,\alpha)} =
[[(\tilde{\alpha},\alpha),(a,x)],(b,y)] .
$$
By Proposition \ref{l}, this is a Lie bracket if and only if
$(\tilde{\alpha},\alpha)$ is a Maurer-Cartan element of
$\g(V,\eta)\o\m$. The homotopy group
$\pi_2(\g(V,\eta),(\tilde\alpha,\alpha))$ is the centre of the Lie algebra
$\bigl(\g(V,\eta)\o\m\bigr)_{(\tilde\alpha,\alpha)}$.

Given $(\tilde{\alpha},\alpha)\in\MC(\g(V,\eta)\o\m)$, the Lie algebra
$\bigl(\g(V,\eta)\o\m\bigr)_{(\tilde{\alpha},\alpha)}$ is naturally
isomorphic to the central extension of the Lie algebra
$(V^\Dual\oplus\C)\o\m$ with bracket
$$
[x,y]_\alpha = [[(0,\alpha),(0,x)],(0,y)]
$$
associated to the 2-cocyle $\tilde{\alpha}$. This proves the following
result.
\begin{theorem} \label{Hydrodynamic}
Let $\m$ be a commutative ring. The Maurer-Cartan elements of $\g\o\m$ are
elements $(\tilde\alpha,\alpha)\in(\Lambda^2V\oplus\Lambda^3V)\o\m$ such
that the bilinear operation $[-,-]_\alpha$ on $V^\Dual\o\m$ defined by is a
Lie bracket, and $\tilde\alpha$ is a 2-cocycle on the Lie algebra
$(V^\Dual\o\m,[-,-]_\alpha)$.
\end{theorem}

The inhomogeneous Euclidean group $\exp(\iso(V,\eta)\o\m)$ is the
semidirect product of the homogenous Euclidean group
$\exp(\so(V,\eta)\o\m)$ and the translation group $V\o\m$. The group
$\exp(\so(V,\eta)\o\m)$ acts on $\MC(\g(V,\eta)\o\m)$ through its adjoint
action on $V^\Dual\o\m$, while $V\o\m$ acts on $\MC(\g(V,\eta)\o\m)$ by
shifting the 2-cocycle $\tilde\alpha$: if $v\in V\o\m$,
$$
v \* (\tilde\alpha,\alpha) = \bigl( \tilde\alpha+v(\alpha) , \alpha \bigr)
,
$$
where $v(\alpha)(x,y)=v([x,y]_\alpha)$. The quotient of
$\MC(\g(V,\eta)\o\m)$ by this group action is $\pi_0(\g(V,\eta)\o\m)$.

The group $\pi_1(\g(V,\eta)\o\m,(\tilde\alpha,\alpha))$ is the quotient of
the subgroup of $\exp(\iso(V,\eta)\o\m)$ consisting of automorphisms of the
Lie algebra $(\g(V,\eta)\o\m)_{(\tilde\alpha,\alpha)}$ by inner ones.

\section{Soloviev's Lie bracket in the formal calculus of variations}

Let $P$ be a Poisson tensor on an affine space $V$. Soloviev
\cite{Soloviev} has constructed a Lie bracket on the infinite jet space of
$V$ which prolongs the Poisson bracket of $V$. In this section, we
generalize Soloviev's construction to Poisson supermanifolds.

The main application we have in mind is to the $1$-symplectic supermanifold
$\Om X$ associated to a manifold $X$, whose Poisson algebra is the Schouten
algebra of $X$. This case is far simpler than the general theory, and we
have taken advantage of this at certain places in our exposition, where the
general theory becomes a little complicated. However, just as in the case
of Poisson manifolds, the general case may be reduced to the case $\Om X$.

\subsection{Higher Euler operators on supermanifolds}

Let $\C^{m|n}$ be a superspace, with coordinates $t^a$. Let $|a|=|t^a|$
equal $0$ or $1$ depending on whether $t^a$ is even or odd. If $U$ is an
open subset of $\C^{m|n}$, let $\CO(U)$ be its (graded) ring of holomorphic
functions. Let $\p_a$ be the derivation $\p/\p t^a:\CO(U)\to\CO(U)$.

Let $\CO_\infty(U)$ be the graded commutative algebra
$$
\CO_\infty(U) = \CO(U)[t^a_k \mid k>0] ,
$$
where $|t^a_k|=|a|$. Let $\p_{k,a}$ be the derivation $\p/\p t^a_k :
\CO_\infty(U)\to\CO_\infty(U)$. We write $t^a_0$ for the generators $t^a$
of $\CO(U)\subset\CO_\infty(U)$, and $\p_{0,a}$ for the derivations $\p_a$.

The algebra $\CO_\infty(U)$ is the space of holomorphic functions on the
supermanifold $J_\infty(U)$ of infinite jets of curves in $U$; such a jet
may be parametrized by the formula
$$
t^a(x) = \sum_{k=0}^\infty \frac{x^k}{k!} t^a_k .
$$
The derivation of $\CO_\infty(U)$ representing differentiation with respect
to $x$ plays a fundamental role: it is given by the formula
$$
\p = \sum_{k=0}^\infty t^a_{k+1} \, \p_{k,a} .
$$

Let $\delta_{k,a}:\CO_\infty(U)\to\CO_\infty(U)$ be the higher Euler
operators of Kruskal et al.\ \cite{KMGZ}
$$
\delta_{k,a} = \sum_{i=0}^\infty (-1)^i \tbinom{k+i}{k} \p^i \, \p_{k+i,a}
,
$$
and let
$$
\delta_k = dt^a\,\delta_{k,a} : \CO_\infty(U) \longrightarrow
(\C^{m|n})^\Dual \o \CO_\infty(U)
$$
be the total higher Euler operators. These are not derivations: indeed,
they are infinite-order differential operators. However, unlike the
derivations $\p_{k,a}$, they have simple transformation properties under
changes of coordinates.

\begin{proposition} \label{transform}
If $f:U\to V$ is a holomorphic map between open subsets of $\C^{m|n}$,
there is a unique homomorphism of algebras
$$
f^* : \CO_\infty(V) \longrightarrow \CO_\infty(U)
$$
which extends the homomorphism $f^*:\CO(V)\to\CO(U)$ and satisfies
$\p\*f^*=f^*\*\p$.

Let $J=df\in\End(\C^{m|n})\o\CO(U)$ be the Jacobian of $f$. For
$u\in\CO_\infty(V)$ and $k\ge0$,
$$
\delta_{k,a} (f^*u) = J_a^b \, f^*(\delta_{k,a}u) \in \CO_\infty(U) .
$$
\end{proposition}
\begin{proof}
It suffices to define $f^*$ on the generators $x_k^a$ of $\CO_\infty(V)$
over $\CO(V)$. By the hypotheses on $f^*$, we have
$$
f^*t_k^a = f^* \p^k t_a = \p^k f^*t_a ,
$$
so the definition of $f^*$ is forced.

By induction on $\ell$, we see that
$$
\p_{k,a} \* \p^\ell = \sum_{j=0}^\ell \tbinom{\ell}{j} \, \p^{\ell-j} \*
\p_{k-j,a} .
$$
It follows that $\p_{k,a} f^*x^b_\ell=\tbinom{\ell}{k} \p^{\ell-k} J_a^b$,
and hence that, for $u\in\CO_\infty(V)$,
$$
\p_{k,a} (f^*u) = \sum_{\ell=k}^\infty \tbinom{\ell}{k}
\bigl(\p^{\ell-k}J^b_a\bigr) \, f^* (\p_{\ell,b}u) .
$$
Thus
\begin{align*}
\delta_{k,a} (f^*u) &= \sum_i (-1)^i \tbinom{k+i}{k} \p^i \bigl(\p_{k+i,a}
f^*u\bigr) \\
&= \sum_{i,\ell} (-1)^i \tbinom{k+i}{k} \tbinom{\ell}{k+i} \p^i \Bigl(
\bigl( \p^{\ell-k-i} J^b_a \bigr) f^*\bigl(\p_{\ell,b}u\bigr) \Bigr) \\
&= \sum_{i,j,\ell} (-1)^i \tbinom{k+i}{k} \tbinom{\ell}{k+i} \tbinom{i}{j}
\bigl( \p^{\ell-k-j} J^b_a \bigr) f^*\bigl(\p^j\p_{\ell,b}u\bigr) \\
&= \sum_{i,j,\ell} (-1)^i \tbinom{\ell-k-j}{i-j} \tbinom{\ell}{k}
\tbinom{\ell-k}{j} \bigl( \p^{\ell-k-j} J^b_a \bigr)
f^*\bigl(\p^j\p_{\ell,b}u\bigr) .
\end{align*}
The sum over $i$ reduces to $\delta(i,j)\,\delta(\ell,j+k)$, and the
right-hand side to $J^b_a (f^*\delta_{k,b}u)$.
\end{proof}

Now suppose that $f$ is a diffeomorphism, and define
$$
f_*=(f^{-1})^*:\CO_\infty(U)\to\CO_\infty(V) .
$$
Since $(gf)_*=g_*f_*$, it follows that $\CO_\infty(U)$ is a module over the
pseudo(super)group of holomorphic diffeomorphisms between open subsets of
$\C^{m|n}$. Thus, the definition of the sheaf of graded commutative
algebras $\CO_\infty$ extends to any $(m|n)$-dimensional complex
supermanifold $M$, and, by Proposition \ref{transform}, the higher Euler
operators extend as well: $\delta_0$ is a connection on the $\CO$-module
$\CO_\infty$, and the higher variational derivatives $\delta_k$, $k>0$, are
sections of $\Om^1\o_\CO\End_\CO(\CO_\infty)$.

\subsection{Soloviev's bracket}

Let $P=P^{ab}\p_a\o\p_b$ be a $\nu$-Poisson tensor on an open subset $U$ of
the superspace $\C^{m|n}$. The following bracket on $\CO_\infty(U)$ was
introduced by Soloviev \cite{Soloviev} (although he restricts attention to
the case $\nu=0$):
\begin{equation} \label{Soloviev}
\{u,v\} = - \sum_{k,\ell} (-1)^{(|b|+\nu)|u|} \p^{k+\ell} \bigl( P^{ab}
\delta_{k,a}u \, \delta_{\ell,b}v \bigr) .
\end{equation}
It is obvious that this bracket extends the Poisson bracket on the subspace
$\CO(U)$ of $\CO_\infty(U)$. However, unlike the Poisson bracket on
$\CO(U)$, Soloviev's bracket does not act by derivations; this is a
fundamental difference between the Hamiltonian formalisms for mechanics and
field theory.

It follows from Proposition \ref{transform} that the bracket
\eqref{Soloviev} is invariant under changes of coordinate; hence the
definition of the Soloviev bracket extends to the sheaf $\CO_\infty$ on a
holomorphic $\nu$-Poisson supermanifold $(M,P)$.

\begin{lemma} \label{bracketd}
$\p\{u,v\} = \{u,\p v\}$
\end{lemma}
\begin{proof}
From the formula
$$
[\p_{k,a},\p] = \begin{cases} \p_{k-1,a} & k>0 , \\ 0 & k=0 ,
\end{cases}$$
it follows that $\delta_{0,a}\p=0$ and that
$\delta_{k,a}\p=\delta_{k-1,a}$ for $k>0$; the lemma follows easily from
this formula.
\end{proof}

Since we are only interested in the case where $M$ is the $1$-symplectic
supermanifold $\Om X$ associated to a manifold $X$, it suffices for our
purposes to extend Soloviev's proof that his bracket satisfies the Jacobi
rule to $\nu$-Poisson tensors $P$ in the special case that their
coefficients $P^{ab}$ are constant. The general case may be reduced to this
one, by expressing the Poisson bracket for a general $\nu$-Poisson tensor
in terms of the Schouten bracket.

The first step in the proof is the following remarkable identity (Statement
6.1.1 of \cite{Soloviev}).
\begin{lemma} \label{constant}
If the coefficients $P^{ab}$ are constant, then
$$
\{u,v\} = - \sum_{k,\ell} (-1)^{(|b|+\nu)|u|} P^{ab} \bigl(
\p^\ell\p_{k,a}u \bigr) \, \bigl( \p^k\p_{\ell,b}v \bigr) .
$$
\end{lemma}
\begin{proof}
We have
\begin{align*}
\{u,v\} &= - \sum_{i,j,k,\ell} (-1)^{i+j+(|b|+\nu)|u|} \tbinom{k+i}{k}
\tbinom{\ell+j}{\ell} P^{ab} \p^{k+\ell} \bigl( \p^i \p_{k+i,a}u \bigr) \,
\bigl( \p^j \p_{\ell+j,b}v \bigr) \\
&= - \sum_{i,j,k,\ell,p} (-1)^{i+j+(|b|+\nu)|u|} \tbinom{k}{k-i}
\tbinom{\ell}{\ell-j} \tbinom{k+\ell-i-j}{k+p-i} P^{ab} \bigl( \p^{k+p}
\p_{k,a}u \bigr) \, \bigl( \p^{\ell-p} \p_{\ell,b}v \bigr) .
\end{align*}
Since $\sum_i (-1)^i \tbinom{k}{k-i} \tbinom{n-i}{m-i} = \tbinom{n-k}{m}$,
this in turn equals
$$
- \sum_{j,k,\ell,p} (-1)^{j+(|b|+\nu)|u|} \tbinom{\ell}{\ell-j}
\tbinom{\ell-j}{k+p} P^{ab} \bigl( \p^{k+p} \p_{k,a}u \bigr) \, \bigl(
\p^{\ell-p} \p_{\ell,b}v \bigr) .
$$
The sum over $j$ reduces to $\delta(\ell,k+p)$, and the lemma follows.
\end{proof}

We now apply the following lemma.
\begin{lemma}
Suppose that
$$
\{\{u,v\},w\} = \alpha(u|v,w) - (-1)^{(|u|+\nu)(|v|+\nu)} \alpha(v|u,w) ,
$$
where $\{u,v\}$ is an operation of degree $\nu$ satisfying
$\{u,v\}=-(-1)^{(|u|+\nu)(|v|+\nu)} \{v,u\}$ and $\alpha(u|v,w)$ is an
operation of degree $0$ satisfying
$$\alpha(u|v,w)=(-1)^{(|v|+\nu)(|w|+\nu)}\alpha(u|w,v).$$ Then $\{u,v\}$ is a
graded Lie bracket.
\end{lemma}

A lengthy calculation shows that when $P^{ab}$ is constant, the hypotheses
of this lemma hold for the bracket \eqref{Soloviev} on $\CO_\infty(U)$,
with $\alpha(u|v,w)$ given by the formula
\begin{multline*}
\alpha(u|v,w) = - \sum_{i,j,k,\ell,p,q}
(-1)^{|b|\,|u|+|d|\,|u|+(|d|+\nu)(|v|+\nu)+(a+b+\nu)c} \\ \tbinom{j}{p}
\tbinom{\ell}{q} P^{ab} P^{cd} \, \bigl( \p^{j+\ell-p-q} \p_{k,c} \p_{i,a}u
\bigr) \, \bigl( \p^{i+q} \p_{j,b}v \bigr) \, \bigl( \p^{k+p} \p_{\ell,d}w
\bigr) .
\end{multline*}

Let $\CP$ be the cokernel of the derivation $\p : \CO_\infty \to
\CO_\infty$, and denote the natural projection from $\CO_\infty$ to $\CP$
by $u\mapsto\tint u \, dx$. Lemma \ref{bracketd} implies that the Lie
bracket $\{u,v\}$ on $\CO_\infty$ induces a graded Lie bracket on $\CP$,
given by the formula
\begin{equation} \label{CP}
\{ \tint u \,dx , \tint v\,dx \} = \tint \{u,v\} \, dx = -
(-1)^{(|b|+\nu)|u|} \tint P^{ab} \delta_au \, \delta_bv \, dx .
\end{equation}

\subsection{The Schouten bracket}
Let $\pi:\Om X\to X$ be the projection from $\Om X$ to $X$, denote the
sheaf $\pi_*\CO_\infty$ on $X$ by $\Lambda_\infty$, and its bracket by
$[u,v]_\infty$. The grading $\Lambda_\infty$ is shifted by $-1$ in the same
way as the grading of $\Lambda$: sections of $\Lambda_\infty^p$ are those
with $p+1$ factors of $\theta_{k,a}$.

In a coordinate system of the form $\{t^a,\theta_a\}$, the Poisson tensor
\eqref{Om} is constant; applying Lemma \ref{constant}, we obtain the
following formula for the bracket on $\Lambda_\infty$:
\begin{equation} \label{[]}
[u,v]_\infty = \sum_{k,\ell} \tint \bigl( (-1)^{|u|} \p^\ell \p^a_ku \*
\p^k\p_{\ell,a}v - \p^\ell\p_{k,a}u \* \p^k\p^a_\ell v \bigr) \, dx .
\end{equation}
Note that the inclusion $\Lambda\hookrightarrow\Lambda_\infty$ is a
morphism of graded Lie algebras.

In the special case where $M$ equals $\Om X$, the bracket \eqref{CP} on
$\CP$ is the Schouten bracket of the formal calculus of variations,
introduced by Gelfand and Dorfman \cite{GD} and Olver \cite{Olver}. Denote
the sheaf $\pi_*\CP$ on $X$ by $\L$, and its bracket by $\[u,v\]$; we grade
$\L$ in the same way as the sheaves $\Lambda$ and $\Lambda_\infty$. As a
graded Lie algebra, $\L$ is a quotient of $\Lambda_\infty$, and
$\tint:\Lambda_\infty\to\L$ is a morphism of graded Lie algebras. The
Schouten bracket is given by two rather different formulas,
\begin{align*}
\[u,v\] &= \tint ( (-1)^{|u|} \delta^au \* \delta_av - \delta_au \*
\delta^av ) \, dx \\
&= \sum_{k,\ell} \tint \bigl( (-1)^{|u|} \p^\ell \p^a_ku \*
\p^k\p_{\ell,a}v - \p^\ell\p_{k,a}u \* \p^k\p^a_\ell v \bigr) \, dx ,
\end{align*}
the first of which manifests the invariance of the bracket under coordinate
transformations, while the second seems to be easier to apply in explicit
calculations.

\section{Hamiltonian manifolds}

The characterization of Hamiltonian operators via the Maurer-Cartan
equation is due to Gelfand and Dorfman \cite{GD}. The following is a global
form of their definition.
\begin{definition}
A \emph{Hamiltonian manifold} $(X,\CQ)$ is a manifold $X$ together with a
section $\CQ\in\Gamma(X,\L^1)$ satisfying the Maurer-Cartan equation
$\[\CQ,\CQ\]=0$. The section $\CQ$ is called a \emph{Hamiltonian
operator}.
\end{definition}

A Hamiltonian operator has a canonical form $\CQ = \tint \theta_a
\DD^{ab} \theta_b \, dx$, where
$$
\DD^{ab} = \sum_{k=0}^N \DD^{ab}_k \, \p^k
$$
is a formally skew-adjoint system of ordinary differential operators with
coefficients in the sheaf $\CO_\infty$. Formal skew-adjointness means that
for every section $u$ of the sheaf $\CO_\infty$,
$$
\sum_{k=0}^N \bigl( \DD^{ab}_k \bigl( \p^ku \bigr) + (-\p)^k \bigl(
\DD^{ba}_ku \bigr) \bigr) = 0 .
$$
For example, if $X=\C$ and $\CQ=\tint\theta(\frac{1}{8}\p^3+t\p)\theta\,dx$
(the second Hamiltonian operator of the KdV hierarchy), the operator $\DD$
equals $\frac{1}{8}\p^3+t\,\p+\half\p t$.

The analogue of Proposition \ref{L} holds for Hamiltonian operators:
$\CQ\in\Gamma(X,\L^1)$ is a Hamiltonian operator if and only if the
morphism of graded sheaves $\delta_\CQ = \[\CQ,-\]$ on $\L$ is a
differential. Denote the sheaf of dg Lie algebras $\L$ with this
differential by $\L_\CQ$; it controls deformations of $\CQ$ in the same way
that the sheaf of dg Lie algebras $\Lambda_Q$ on a Poisson manifold
controls deformations of the Poisson tensor $Q$.
\begin{definition}
A Hamiltonian operator $\CQ$ is \emph{Darboux} if the sheaf of dg Lie
algebras $\L_\CQ$ is formal.
\end{definition}

\subsection{A resolution of $\L$}

We now introduce a resolution $\LL$ of the sheaf of graded Lie algebras
$\L$; this resolution is a sheaf of Fock spaces.

Let $\tilde\Lambda_\infty=\Lambda_\infty/(\C\*1)$ be the quotient of
$\Lambda_\infty$ by its centre, and let $\LL$ be the cone of the morphism
$\p:\tilde\Lambda_\infty\to\Lambda_\infty$; in other words, $\LL$ is
isomorphic to the graded sheaf
$\Lambda_\infty\oplus\tilde\Lambda_\infty[1]$, where $\tilde\CO_\infty[1]$
is a copy of $\tilde\Lambda_\infty$ shifted in degree by $-1$. Denoting
elements of $\tilde{\Lambda}_\infty[1]$ by $\eps\tilde{u}$, the
differential equals $D(u+\eps\tilde{u}) = \p\tilde{u}$. Equipped with the
bracket
$$
[ u+\eps\tilde{u} , v+\eps\tilde{v} ]_\infty = [ u,v ]_\infty + \eps \bigl(
[ \tilde{u},v ]_\infty + (-1)^{|u|} [ u,\tilde{v} ]_\infty \bigr) ,
$$
$\LL$ is a sheaf of dg Lie algebras.

\begin{theorem} \label{LL}
The morphism $\tint:\LL\to\L$ defined by the formula
$$
\tint (u+\eps\tilde{u}) = \tint u \, dx ,
$$
is a weak equivalence of dg Lie algebras.
\end{theorem}
\begin{proof}
It is clear that $\tint$ is compatible with the differential on $\LL$:
$$
\tint D (u+\eps\tilde{u}) = \tint (\p \tilde{u}) \, dx = 0 .
$$
It is also easy to see that $\tint$ is a morphism of graded Lie algebras,
since
$$
\tint [ u+\eps\tilde{u} , v+\eps\tilde{v} ]_\infty = \[
\tint(u+\eps\tilde{u}) , \tint(v+\eps\tilde{v}) \] .
$$
It only remains to check that $\tint$ is a weak equivalence; this is a
variant on the ``exactness of the variational bicomplex.'' We learned the
idea used in the following proof from E.~Frenkel.

Let $U$ be a connected open subset of $\C^{m|n}$, and let
$u\in\CO_\infty(U)$. We must show that $\p u=0$ if and only if $u$ is a
multiple of $1$. It is clear that this is so if $u\in\CO(U)$, since in that
case, $\p u = \p t^a \, \p_au$. The operators $\p$ and
$$
\rho = \sum_{k=0}^\infty k(k+1) \, t_k^a \, \p_{k+1,a}
$$
generate an action of the Lie algebra $\Sl(2)$ on $\CO_\infty(U)$, whose
Cartan subalgebra acts by the semisimple endomorphism
$$
H = \sum_{k=0}^\infty k \, t^a_k \, \p_{k,a} ,
$$
with kernel $\CO(U)$.

Suppose that $\p u=0$. Since $\rho^iu=0$ for $i\gg0$, we see that the
irreducible $\Sl(2)$-module spanned by $u$ is finite-dimensional. Since a
finite-dimensional representation of $\Sl(2)$ on which $H$ has non-negative
spectrum is trivial, we conclude that $Hu=0$; hence $u$ lies in
$\CO(U)\subset\CO_\infty(U)$, and as we have seen, is a multiple of $1$.
\end{proof}

\subsection{Ultralocal Hamiltonian operators}
Since the bracket on $\L$ prolongs the Schouten bracket, a Poisson tensor
$Q$ on $X$ gives rise to a Hamiltonian operator $\CQ=\tint Q\,dx$.  Such
Hamiltonian operators are called \emph{ultralocal}.
\begin{theorem} \label{ultra}
If $Q$ is the Poisson tensor associated to a symplectic manifold $(X,\om)$,
the inclusion of sheaves of dg Lie algebras
$\Lambda_Q\hookrightarrow\L_\CQ$ is a weak equivalence; in particular, the
Hamiltonian operator $\CQ=\tint Q\,dx$ is Darboux.
\end{theorem}
\begin{proof}
We must show that if $(X,Q)$ is a symplectic manifold, the inclusion
$$
(\Lambda^\bull,\delta_Q)\hookrightarrow(\LL^\bull,D+\delta_Q)
$$
of sheaves of dg Lie algebras is a weak equivalence. By the Darboux theorem
(in its original sense!), it suffices to consider a convex subset $U$ of
$\C^{2\ell}$ with its standard symplectic structure, and Poisson tensor
$$
Q = \sum_{a=1}^\ell \theta_a \, \theta_{a+\ell} .
$$
Let $\delta_\Q$ be the differential associated to the Maurer-Cartan element
$$
\Q = \sum_{a=1}^\ell \theta_{0,a} \, \theta_{0,a+\ell} .
$$
of $\LL(U)$; it is given by the formula
$$
\delta_\Q = \sum_{k=0}^\infty \sum_{a=1}^\ell \bigl( \theta_{k,a}
\, \p_{k,a+\ell} - \theta_{k,a+\ell} \, \p_{k,a} \bigr) .
$$
Clearly, the dg Lie algebra $\LL_\Q(U)=(\LL(U),D+\delta_\Q)$ is a
resolution of $(\L,\delta_{\int\!Q\,dx})$.

The complex $\LL_\Q(U)$ is isomorphic to the cone of the morphism
$$
\p : \tilde{\Om}^\bull(J_\infty(U)) \longrightarrow \Om^\bull(J_\infty(U))
,
$$
where $\Om^\bull(J_\infty(U))$ is the de Rham complex of the jet-space
$J_\infty(U)$ and $\tilde{\Om}^\bull(J_\infty(U))$ is its quotient by the
constant functions. To see this, one identifies $\theta_{k,a}$ with
$dt_k^{a+\ell}$ and $\theta_{k,a+\ell}$ with $-dt_k^a$. Theorem \ref{ultra}
now follows from the de Rham theorem for $J_\infty(U)$.
\end{proof}

If the Poisson tensor $Q$ is not symplectic, the inclusion
$\Lambda_Q\hookrightarrow\L_\CQ$ is not a weak equivalence; this is obvious
if the Poisson tensor vanishes, and the general case may be inferred from
this one.

\subsection{Hamiltonian manifolds of hydrodynamic type}

Let $\eta$ be a flat contravariant metric on $M$, with coefficients
$\eta^{ab}=\eta(dt^a,dt^b)$. Dubrovin and Novikov \cite{DN} associate to
$\eta$ a Hamiltonian operator $H_\eta$; in flat coordinates (those for
which the coefficients $\eta^{ab}$ are constant), it is given by the
formula
$$
H_\eta = \half \tint \eta^{ab} \theta_a \, \p\theta_b \, dx .
$$
The differential $d_\eta=\[H_\eta,-\]$ on $\L$ is given by the formula
$$
d_\eta \tint u \, dx = - \sum_k \eta^{ab} \tint \theta_{k+1,a} \p_{k,b} u
\, dx ,
$$
and the resulting sheaf of dg Lie algebras is denoted $\L_\eta$. We may now
state the main result of this paper. Let $\g(X,\eta)$ be the sheaf of
graded Lie algebras on $X$ whose stalk at $x\in X$ is the graded Lie
algebra $\g(T^\Dual_xX,\eta)$ introduced in Section
\ref{hydro-example}. Let $\tau_0:\g(U,\eta)_x\to\Lambda_\infty(U)$ be the
operation which substitutes $\theta^a_0$ for $\theta^a$.
\begin{theorem} \label{hydro}
The morphism $\sigma:\g(X,\eta)\hookrightarrow\L_\eta$ defined by the
formula
$$
\sigma ( \tilde\alpha , \alpha ) = \tint \tau_0(\tilde\alpha) \, dx + \tint
\eta_{ab} \, t_0^a \, \p_0^b \tau_0(\alpha) \, dx
$$
is a weak equivalence of sheaves of dg Lie algebras. In particular,
hydrodynamic Hamiltonian operators are Darboux, and $\sigma$ induces a weak
equivalence of sheaves of Deligne 2-groupoids
$$
\CC(\sigma) : \CC(\g(X,\eta)) \simeq \CC(\L_\eta) .
$$
\end{theorem}

\subsection{Lifting Hamiltonian operators to $\LL$}
The proof of Theorem \ref{ultra} was based on the idea of lifting the
Hamiltonian operator $\CQ=\tint Q\,dx$ to a Maurer-Cartan element of
$\LL$. This may be generalized as follows.
\begin{definition}
A \emph{lift} of a Hamiltonian manifold $(X,\CQ)$ is a section $\Q$ of
$\LL^1$ with $\CQ=\tint\Q\,dx$, and which satisfies the Maurer-Cartan
equation $D\Q+\half[\Q,\Q]_\infty=0$.
\end{definition}

If $\Q$ is a lift of a Hamiltonian operator $\CQ$, there is a weak
equivalence of sheaves of dg Lie algebras
$$
\tint : \bigl( \LL , D + \delta_\Q \bigr) \longrightarrow \bigl( \L ,
\delta_\CQ \bigr) ,
$$
where $\delta_\Q$ is the differential $\delta_\Q u=[\Q,u]_\infty$ on
$\LL$.

Let us give some explicit examples of lifts. As we have already observed,
an ultralocal Hamiltonian operator $\tint Q\,dx$ has the lift $Q$. A
hydrodynamic Hamiltonian operator
$\half\tint\eta^{ab}\theta_a\p\theta_b\,dx$ has the lift
$\half\eta^{ab}\theta_a\theta_{1,b}$. Since a manifold $X$ with flat
contravariant metric $\eta$ has an atlas whose charts are flat and whose
transition functions are inhomogeneous orthogonal transformations, these
lifts patch together to give a lift of $H_\eta$ over all of $X$.

For a less trivial example, the second Hamiltonian operator of the KdV
hierarchy, $\CQ = \tint \theta \bigl( \frac18\p^3 + t \, \p \bigr) \theta
\, dx$, has a family of lifts (cf.\ Dickey \cite{Dickey})
$$
\Q = \tfrac18\theta\theta_3 + t\theta\theta_1 + a \, \p(\theta\theta_2) +
\tfrac18\eps\,\theta\theta_1\theta_2 , \quad a \in \C .
$$
\begin{proposition}
Every Hamiltonian manifold $(X,\CQ)$ which is Stein has a lift $\Q$.
\end{proposition}
\begin{proof}
Lifts $\Q=u+\eps\tilde{u}$ of $\CQ$ are characterized by the equations
$\CQ=\tint u\,dx$ and
$$
\p\tilde{u} + \half[u,u]_\infty = [u,\tilde{u}]_\infty = 0 .
$$
Let $u$ be a section of $\Lambda^1_\infty$ such that $\CQ=\tint u\,dx$;
there are no obstructions to the existence of $u$, because $X$ is
Stein. Since $\[\CQ,\CQ\]=\tint[u,u]_\infty\,dx=0$, we see that there is a
section $\tilde{u}$ of $\Lambda^2_\infty$ such that
$\p\tilde{u}+\half[u,u]_\infty=0$; again, there are no obstructions to the
existence of $\tilde{u}$. Taking the bracket of this equation with $u$, we
see that
$$
[u,\p\tilde{u}]_\infty+ \half [u,[u,u]_\infty]_\infty = 0 .
$$
But $[u,[u,u]_\infty]_\infty$ vanishes by the Jacobi rule, while
$[u,\p\tilde{u}]_\infty=\p[u,\tilde{u}]_\infty$. By Theorem \ref{LL}, we
conclude that $[u,\tilde{u}]_\infty=0$.
\end{proof}

\section{The proof of the main theorem}

We now give the proof of Theorem \ref{hydro}. The hydrodynamic Hamiltonian
operator $\half\tint\eta^{ab}\theta_a\p\theta_b\,dx$ has the lift
$\half\eta^{ab}\theta_a\theta_{1,b}$. The associated differential
of $\LL(U)$ equals
$$
d_\eta = \[\half\eta^{ab}\theta_a\theta_{1,b},-\] = -d+\half\p\*d_0 ,
$$
where $d = \sum_{k=0}^\infty \eta^{ab} \theta_{k+1,a} \p_{k,b}$ and $d_0 =
\sum_{k=0}^\infty \eta^{ab} \theta_{k,a} \p_{k,b}$.
\begin{lemma} \label{tau}
Let $\eta$ be a constant metric on $\C^n$, and let $U$ be a convex subset
of $\C^n$ containing $0$. The map of graded vector spaces
$\tau:\g(U,\eta)_0=\g(T_0U,\eta_0)\to\LL_\eta(U)$, defined on
$\g^p(U,\eta)_0$ by the formula
$$
\tau ( \tilde\alpha , \alpha ) = \tau_0(\tilde\alpha) + ( \eta_{ab} \,
t_0^a \, \p_0^b - \half \eps \, p ) \tau_0(\alpha) ,
$$
is a morphism of dg Lie algebras.
\end{lemma}
\begin{proof}
1) $\tau$ is a morphism of complexes (that is, $(D+d_\eta)\*\tau=0$): Let
$(\tilde\alpha,\alpha)$ be an element of $\g^p(U,\eta)_0$. It is obvious
that
$$
(D+d_\eta)\tau(\tilde\alpha,0) = (D+d_\eta)\tau_0(\tilde\alpha) = 0 ,
$$
since $D\tau_0(\tilde\alpha)$ and $d_\eta\tau_0(\tilde\alpha)$ both
vanish. As for $(D+d_\eta) \tau ( 0 ,\alpha )$, we have
$$
D \bigl( \eta_{ab} \, t_0^a \, \p_0^b \bigr) \tilde\alpha = d_\eta \bigl( -
\half \, \eps \, p \bigr) \tilde\alpha = 0
$$
and
$$
d_\eta \bigl( \eta_{ab} \, t_0^a \, \p_0^b \bigr) \tilde\alpha + D \bigl( -
\half \, \eps \, p \bigr) \tilde\alpha = \half \, p \, \p \, \tilde\alpha -
\half \, p \, \p \, \tilde\alpha = 0 .
$$

2) $\tau$ preserves the Lie bracket: If $\alpha\in\h^p(U,\eta)_0$, we have
\begin{multline*}
\tau [ (0,\alpha) , (\tilde\beta,\beta) ] = \tau( - \eta_{ab} \,
\p^a\alpha \, \p^b\tilde\beta , (-1)^{p+1} \eta_{ab} \, \p^a\alpha \,
\p^b\beta ) \\
\begin{aligned}
{}&= - \tau_0( \eta_{ab} \, \p^a\alpha \, \p^b\tilde\beta ) +
(-1)^{p+1} \, \eta_{ab} \bigl( \eta_{cd} \,  t^c_0 \, \p^d_0 - \half \,
\eps \, (p+q) \bigr) \tau_0 \bigl( \p^a\alpha \, \p^b\beta \bigr) \\
&= - \eta_{ab} \, \p_0^a\tau_0(\alpha) \, \p_0^b\tau_0(\tilde\beta) +
(-1)^{p+1} \, \eta_{ab} \bigl( \eta_{cd} \, t^c_0 \, 
\p^d_0 - \half \, \eps \, (p+q) \bigr) \p_0^a\tau_0(\alpha) \,
\p_0^b\tau_0(\beta) .
\end{aligned}
\end{multline*}
On the other hand,
\begin{align*}
[ \tau(0,\alpha) , \tau(\tilde\beta,\beta) ]_\infty &= [ ( \eta_{ab} \,
t_0^a \, \p_0^b - \half\eps \, p ) \tau_0(\alpha) , \tau_0(\tilde\beta) + (
\eta_{cd} \, t_0^c \, \p_0^d - \half\eps \, q ) \tau_0(\beta) ]_\infty \\
&= [ \eta_{ab} \, t_0^a \, \p_0^b \tau_0(\alpha) , \tau_0(\tilde\beta) +
\eta_{cd} \, t_0^c \, \p_0^d \tau_0(\beta) ]_\infty \\
&- \half \, \eps \bigl( (-1)^p \, q \, [ \eta_{ab} \, t_0^a \, \p_0^b
\tau_0(\alpha) , \tau_0(\beta) ]_\infty +  p \, [ \tau_0(\alpha) ,
\eta_{cd} \, t_0^c \, \p_0^d \tau_0(\beta) ]_\infty \bigr) \\
&= - \eta_{ab} \, \p^b_0\tau_0(\alpha) \, \p^a_0\tau_0(\tilde\beta) \\
&+ (-1)^p \eta_{ab} \eta_{cd} \, t_0^a \, \p_0^c \, \p_0^b \tau_0(\alpha)
\, \p_0^d \tau_0(\beta)
- \eta_{ab} \eta_{cd} \, t_0^c \, \p_0^b \tau_0(\alpha) \, \p_0^a \,
\p_0^d \tau_0(\beta) \\
&+ \half \, \eps \bigl( (-1)^p \, q \, \eta_{ab} \,
\p_0^b \tau_0(\alpha) \, \p^a_0\tau_0(\beta)
- (-1)^{p+1} p \, \eta_{cd} \, \p^c_0\tau_0(\alpha) \,
\p_0^d \tau_0(\beta) \bigr) .
\end{align*}
From these formulas, we see that $\tau [ (0,\alpha),(\tilde\beta,\beta)] =
[\tau(0,\alpha),\tau(\tilde\beta,\beta)]_\infty$. Finally, it is clear that
$[\tau(\tilde\alpha,0),\tau(\tilde\beta,0)]_\infty=0$, as they must, since
$(\tilde\alpha,0)$ and $(\tilde\beta,0)$ commute in $\g(U,\eta)_0$.
\end{proof}

The operations $\iota(u+\eps\tilde{u})=\tilde{u}$ and
$\eps(u+\eps\tilde{u})=\eps u$ on $\LL(U)$ satisfy the canonical graded
commutation relations $[\eps,\iota]=1$; using $\iota$, the differential of
$\LL$ may be written $D=\iota\p$.

\begin{lemma}
The morphism $\RR=1+\half\eps d_0$ (of complexes, not of dg Lie algebras)
induces an isomorphism of complexes
$$
\RR : \LL_\eta(U) \longrightarrow (\LL(U),D-d) .
$$
\end{lemma}
\begin{proof}
We must show that $\RR(D+d_\eta)-(D-d)\RR$ vanishes. Rewritten using the
operators $\eps$ and $\iota$, and taking into account that the operators
$d$ and $d_0$ graded commute with $\iota$ and $\eps$, and that
$[\p,d_0]=0$, we see that this equals
$$
(1+\half\eps d_0)(\iota\p-d+\half\p d_0) - (\iota\p-d)(1+\half\eps d_0) =
\half\p \bigl( d_0 - [\eps,\iota] d_0 + \tfrac14 \eps [d_0,d_0] \bigr) -
\half \eps [d_0,d] ,
$$
which vanishes, since $[d,d_0]=[d_0,d_0]=0$ and $[\eps,\iota]=1$.
\end{proof}

\begin{lemma}
The morphism of complexes $\RR\*\tau:\g(U,\eta)_0\to(\LL(U),D-d)$ has the
formula $\RR\*\tau ( \tilde\alpha , \alpha ) = \tau_0(\tilde\alpha) + \eps
\, \tau_0(\alpha)$.
\end{lemma}
\begin{proof}
We have
\begin{align*}
\RR\*\tau ( \tilde\alpha , \alpha ) &= (1+\half\,\eps\,d_0)\*\bigl(
\tau_0(\tilde\alpha) + \eta_{ab} \, t_0^a \, \p_0^b \tau_0(\alpha) - \half
\eps \, p \, \tau_0(\alpha) \bigr) \\ &= \tau_0(\tilde\alpha) + \eta_{ab}
\, t_0^a \, \p_0^b \tau_0(\alpha) + \half \, \eps \, d_0 \bigl( \eta_{ab}
\, t_0^a \, \p_0^b \tau_0(\alpha) \bigr) - \half \eps \, p \,
\tau_0(\alpha) .
\end{align*}
The formula follows, since $d_0 \bigl( \eta_{ab} \, t_0^a \, \p_0^b
\tau_0(\alpha)=(p+2)\tau_0(\alpha)$.
\end{proof}

Theorem \ref{hydro} is now a consequence of the following lemma.
\begin{lemma}
The morphism $\RR\*\tau:\g(U,\eta)_0\to(\LL(U),D-d)$ is a weak equivalence.
\end{lemma}
\begin{proof}
There is short exact sequence of complexes
$$
0 \longrightarrow \bigl( \Lambda_\infty(U) , -d \bigr) \longrightarrow
(\LL(U),D-d) \longrightarrow \bigl( \tilde{\Lambda}_\infty(U)[1] , -d
\bigr) \longrightarrow 0 ,
$$
and hence, for $p\ge-1$, a long exact sequence
$$
\cdots \longrightarrow H^{p-1}(\LL(U),D-d) \longrightarrow H^p\bigl(
\tilde{\Lambda}_\infty(U) , d \bigr) \xrightarrow{\ \delta\ } H^p\bigl(
\Lambda_\infty(U) , d \bigr) \longrightarrow \cdots
$$

There is an isomorphism between the complex $(\Lambda_\infty(U),d)$ and the
de~Rham complex $\Om^\bull(J_\infty(U),\Lambda\C^n)[1]$, obtained by
mapping $\theta_{k+1,a}$ to $\eta_{ab}\,dt^b_k$ and $\theta_{0,a}$ to the
basis vector $\theta_a$ of $\C^n$. Likewise, the complex
$(\tilde{\Lambda}_\infty(U),d)$ is isomorphic to the reduced de~Rham
complex $\tilde{\Om}^\bull(J_\infty(U),\Lambda\C^n)[1]$.

The Poincar\'e lemma for $J_\infty(U)$ shows that $\RR\*\tau$ induces
isomorphisms between the groups $H^p\bigl( \Lambda_\infty(U) , d \bigr)$
and $H^p\bigl( \tilde{\Lambda}_\infty(U) , d \bigr)$ and the group
$\Lambda^{p+1}\C^n$. The composition of $\RR\*\tau$ with the boundary map
$\delta : H^p\bigl( \tilde{\Lambda}_\infty(U) , d \bigr) \to H^p\bigl(
\Lambda_\infty(U) , d \bigr)$ vanishes: if $\alpha\in\Lambda^{p+1}\C^n$, we
have
$$
\delta\*\RR\*\tau(\alpha) = D (\eps\,\tau_0(\alpha)) = \p \tau_0(\alpha) =
0 .
$$
We conclude that there is a short exact sequence
$$
0 \longrightarrow \Lambda^{p+1}\C^n \longrightarrow H^p(\LL(U),D-d)
\longrightarrow \Lambda^{p+2}\C^n \longrightarrow 0 ,
$$
and hence that $\RR\*\tau$ is indeed an isomorphism onto the cohomology of
$D-d$.
\end{proof}

\end{document}